\documentstyle[11pt]{article}

\input amssym.def
\input amssym

\def \bs{\backslash}
\def \C{{\Bbb C}}
\def \CF{{\cal F}}
\def \CH{{\cal H}}
\def \CS{{\cal S}}
\def \det{{\rm det}}
\def \dist{{\rm dist}}
\def \Dom{{\rm Dom}}
\def \Eig{{\rm Eig}}
\def \ga{\gamma}
\def \Ga{\Gamma}

\def \jtoinfty{\begin{array}{c} j\rightarrow \infty \\ \longrightarrow	\\ {}
		\end{array}}
\def \la{\lambda}
\def \lap{\triangle}
\def \mod{{\rm mod}}
\def \N{{\Bbb N}}
\def \ph{{\varphi}}
\def \prf{{\bf Proof: }}
\def \qed{\hfill $\Box$

$ $

}
\def \Q{{\Bbb Q}}
\def \ra{\rightarrow}
\def \R{{\Bbb R}}
\def \Re{{\rm Re}}
\def \S{{\cal S}}
\def \tr{{\rm tr}}
\def \vol{{\rm vol}}
\def \Z{{\Bbb Z}}

\begin{document}

\title{Regularized and $L^2$-Determinants} 

\author{Anton Deitmar\\ {\small Math. Inst. d. Univ., INF 288, 69126 Heidelberg, Germany}}
\date{}
\maketitle

\pagestyle{myheadings}
\markright{REGULARIZED AND $L^2$-DETERMINANTS}

\tableofcontents

\newcommand{\rez}[1]{\frac{1}{#1}}
\newcommand{\der}[1]{\frac{\partial}{\partial #1}}
\newcommand{\binom}[2]{\left( \begin{array}{c}#1\\#2\end{array}\right)}

\newcounter{lemma}
\newcounter{corollary}
\newcounter{proposition}
\newcounter{theorem}

\renewcommand{\subsection}{\stepcounter{subsection}\stepcounter{lemma} 
	\stepcounter{corollary} \stepcounter{proposition}
	\stepcounter{conjecture}\stepcounter{theorem}
	\vspace{4pt}
	{\bf \arabic{section}.\arabic{subsection}\hspace{4pt}}}

\newtheorem{conjecture}{\stepcounter{lemma} \stepcounter{corollary} 	
	\stepcounter{proposition}\stepcounter{theorem}
	\stepcounter{subsection}Conjecture}[section]
\newtheorem{lemma}{\stepcounter{conjecture}\stepcounter{corollary}	
	\stepcounter{proposition}\stepcounter{theorem}
	\stepcounter{subsection}Lemma}[section]
\newtheorem{corollary}{\stepcounter{conjecture}\stepcounter{lemma}
	\stepcounter{proposition}\stepcounter{theorem}
	\stepcounter{subsection}Corollary}[section]
\newtheorem{proposition}{\stepcounter{conjecture}\stepcounter{lemma}
	\stepcounter{corollary}\stepcounter{theorem}
	\stepcounter{subsection}Proposition}[section]
\newtheorem{theorem}{\stepcounter{conjecture} \stepcounter{lemma}
	\stepcounter{corollary}\stepcounter{proposition}		
	\stepcounter{subsection}Theorem}[section]
	
$$ $$

\begin{center}
{\bf Introduction}
\end{center}
In the early seventies D. Ray and I. Singer \cite{RS-Rtors} introduced the notion of zeta-regularized determinants.
They used it to define the analytic version of Reidemeister torsion as an alternating product of determinants.
One way to understand analytic torsion is to consider it as a "multiplicative index" of an elliptic complex.
By the $L^2$-index theorem of M. Atiyah \cite{At} this analogy suggests that one should define "$L^2$-torsion" via the use of von Neumann traces and one should ask whether a "multiplicative $L^2$-index theorem" holds, which would say that analytic torsion should coincide with the $L^2$-torsion.
This, however, fails to hold. 
To measure the failure one considers the quotient $\frac{{\rm analytic\ torsion}}{L^2-{\rm torsion}}$. For this number there sometimes is a Lefschetz theorem expressing it as (regularized) sum of local geometric contributions (see \cite{D-holtors}).

Since the torsion is an alternating product of determinants one should more generally consider the quotient $\frac{\rm regularized\ determinant}{L^2-\rm determinant}$. The study of the latter is the aim of the present paper.
We show that in towers the regularized determinant converges to the $L^2$-determinant and that the regularized determinant can generally be expressed as an "Euler-product" of $L^2$-determinants.

The first section is devoted to the study of regularized determinants.
We show an asymptotic formula which is useful to determine the exact shape of a geometric zeta function \cite{D-det}.
In the second section we study $L^2$-determinants and in the third we prove convergence theorems that show that in a tower of coverings converging to the universal covering the regularized determinant converges to the $L^2$-determinant.
Since $L^2$-determinants or $L^2$-torsion numbers are easier to compute than the regularized determinant or the analytic torsion this result shows the generic nontriviality of the regularized determinant of the analytic torsion.
The last section gives an "Euler-product" formula which says that the regularized determinant is an infinite product of equivariant $L^2$-determinants on certain, generally infinite coverings.

\section{Regularized determinants}
\subsection
Let $a=(a_n)_{n\in \N}$ be a sequence of nonzero complex numbers such that all members $a_n$ lie in the complex angle $W_\epsilon := \{ -\pi +\epsilon <\arg z<\pi -\epsilon \}$ for some $\epsilon >0$.
The sequence $a$ is called {\bf admissible} if
\begin{itemize}
\item[a)]
for some natural number $k$ the sum $\sum_n a_n^{-k}$ converges absolutely and
\item[b)]
the theta series
$$
\Theta_{a}(t) := \sum_n e^{-ta_n},
$$ which is convergent for all $t>0$, has an asymptotic expansion
$$
\Theta_{a}(t) \sim \sum_{\nu =0}^\infty c_\nu t^{\alpha_\nu},
$$
as $t\downarrow 0$, where $\alpha_\nu \in \R$, $\alpha_\nu \ra  +\infty$.
\end{itemize}

\subsection {\bf Remarks:} 
- The asymptotic expansion above means that, assumed the $\alpha_\nu$ are ordered increasingly, for any natural number $N$ there is a constant $C$ such that
$$
| \Theta_{a}(t) - \sum_{\nu =0}^N c_\nu t^{\alpha_\nu}|\leq C t^{\alpha_{N+1}}
$$
for all $t\in ]0,1[$.

- Let $(a_n)_{n\in\N}$ be an admissible sequence then the condition $a_n \in W_\epsilon$ and condition a) imply $\Re (a_n)\ra \infty$ for $n\ra \infty$.

- The same two conditions imply that the theta series $\Theta_{a}(t)$ converges absolutely for all $t>0$.

- Instances of admissible sequences are semilattices, i.e. $a_n = a_0+nz$, for some $a_0,z\in W_\epsilon$ or the sequences of eigenvalues of elliptic differential operators (Lemma 1.7.4 in \cite{Gilk}).

\begin{lemma} \label{zetafunct}
Let $(a_n)_{n\in\N}$ be an admissible sequence, then the zeta series
$$
\zeta_{a}(s) := \sum_n a_n^{-s},
$$
where the power is defined by the principal branch of the logarithm, converges in the half plane $\Re (s) >k$ and extends to a meromorphic function on $\C$.
The only poles of $\zeta_{a}(s)$ are simple poles at $s=-\alpha_\nu$ whenever $\alpha_\nu$ is not an integer $\geq 0$.
The residue at $s=-\alpha_\nu$ is $\frac{c_\nu}{\Ga (-\alpha_\nu)}$.
\end{lemma}

\prf
To show the convergence write $a_n$ as $|a_n|e^{i\theta_n}$ and $s=\sigma +it$ with $\sigma , t\in \R$, then
$|a_n^{-s}|=e^{-\Re (s\log a_n)}=e^{-(\sigma \log |a_n| -t\theta_n)}\leq |a_n|^{-\sigma} e^{|t|(\pi -\epsilon)}$, which gives the claim.
To proof the assertion on extension we may write
$$
\zeta_{a}(s) = \frac{1}{\Ga (s)} \int_0^\infty t^{s-1}\Theta_{a}(t) dt
$$
for $\Re (s)$ sufficiently large.
We split this integral into two parts as $\int_0^\infty = \int_0^1 +\int_1^\infty$.
The second part, i.e. $\frac{1}{\Ga(s)} \int_1^\infty t^{s-1} \Theta_{a}(t) dt$ is convergent for all $s$ and thus represents a holomorphic function in $s$.
For the first part let $N$ be any natural number and write
\begin{eqnarray*}
\frac{1}{\Ga(s)}\int_0^1 t^{s-1}\Theta_{a}(t) dt &=& \frac{1}{\Ga(s)}\sum_{\nu =0}^N c_\nu \frac{1}{s+\alpha_\nu} \\
&+& \frac{1}{\Ga(s)}\int_0^1 t^{s-1}\left( \Theta_{a}(t)-\sum_{\nu =0}^N c_\nu t^{\alpha_\nu}\right) dt.
\end{eqnarray*}
The latter integral converges locally uniformly for $\Re (s) >-\alpha_{N+1}$.
Since $\alpha_N \ra \infty$ as $N$ tends to infinity the claim follows.
\qed

\subsection
We want to use the zeta function $\zeta_{a}$ to define the regularized product of the $a_n$.
To motivate this, consider finite products for the moment.
For a finite tupel $a=(a_1,\dots ,a_N)$ we also may consider the zeta function $\zeta_a(s) := \sum_{n=1}^N a_n^{-s}$, which gives an entire function.
We have $\zeta_a'(s) = \sum_{n=1}^N \log (a_n) a_n^{-s}$ and thus
$\prod_{n=1}^N a_n = \exp (-\zeta_a'(0))$.
So we define:

The {\bf regularized product} of the admissible sequence $a=(a_n)_{n\in\N}$ is defined to be the number:
$$
\hat{\prod_n} a_n := \exp (-\zeta_{a}'(0)).
$$

{\bf Example:}
For the sequence $a_n :=n$ we get $\zeta_{a}=\zeta$, the Riemann zeta function.
Because of $\zeta'(0)=-\frac{1}{2}\log (2\pi)$ we get
$$
\hat{\prod_n} n = \sqrt{2\pi}.
$$

We now give some elementary rules concerning regularized products.

\begin{proposition}
\begin{itemize}

\item[i)]
If $\N = A\cup B$ is a disjoint decomposition of the set of natural numbers and each of the sequences $(a_n)_{n\in A}$ and $(a_n)_{n\in B}$ is either finite or admissible then $(a_n)_{n\in\N}$ is admissible and we have
$$
\hat{\prod_{n\in \N}} a_n = \hat{\prod_{n\in A}} a_n \hat{\prod_{n\in B}} a_n.
$$

\item[ii)]
If $a=(a_n)_{n\in\N}$ is admissible and $c>0$ then $(ca_n)_{n\in\N}$ is admissible and
$$
\hat{\prod_n} ca_n = c^{\zeta_{a}(0)} \hat{\prod_n}a_n.
$$

\item[iii)]
If $(a_n)_{n\in\N}$ is admissible and $\Re (s) >0$ is such that $a_n^s$ still lies in some $W_\epsilon$ then $(a_n^s)$ is admissible and
$$
\hat{\prod_n}a_n^s = \left( \hat{\prod_n}a_n\right)^s.
$$
\end{itemize}
\end{proposition}

\prf
i) is clear since the theta series behaves additively.

ii) We have $\Theta_{(ca_n)_{n\in\N}}(t) = \Theta_{(a_n)_{n\in\N}}(ct)$ and therefore
$$
\zeta_{(ca_n)}(s) = \frac{1}{\Ga(s)} \int_0^\infty t^{s-1} \Theta_{(a_n)_{n\in\N}}(ct) dt = c^{-s}\zeta_{(a_n)_{n\in\N}}(s),
$$
so that $\zeta_{(ca_n)_{n\in\N}}'(s) = -(\log c) c^{-s} \zeta_{(a_n)_{n\in\N}}(s) + c^{-s}\zeta_{(a_n)_{n\in\N}}'(s)$, which gives the claim.

iii) This follows from $\zeta_{(a_n^s)_{n\in\N}}(s') = \zeta_{(a_n)_{n\in\N}}(ss')$.
\qed

\subsection
Now suppose given a Hilbert space $\CH$ and a densely defined operator $A$ on $\CH$.
We call $A$ {\bf regular} if there is a dense subspace $\CH'$ of ${\rm dom} A\subset\CH$ such that $\CH'=\bigoplus_i H_i$, where the $H_i$ are pairwise orthogonal, each $H_i$ is finite dimensional and stable under $A$. 
Further, if $f\in \Dom (A)$ with $f=\sum_i f_i,\ f_i\in H_i$, then we insist that $Af =\sum_i Af_i$.

It is easy to see that regular operators are closable.
Let $(a_n)_{n\in \N}$ be the sequence of eigenvalues of $A$ on $\CH'$, each occurring with its algebraic multiplicity.
If the sequence $(a_n)$ is admissible we call $A$ {\bf admissible} and define its {\bf determinant} to be
$$
\det (A) := \hat{\prod_n}a_n.
$$

{\bf Remarks:}
- Since an admissible sequence tends to infinity it actually follows that for an admissible operator $A$ the space $\CH'$ coincides with the space of $A$-finite vectors, hence is uniquely determined.

- The Proposition 1.3 immediately translates to the rules
\begin{itemize}
\item
$\det (A\oplus B)=\det (A) \det (B)$
\item
$\det(cA) = c^{\zeta (0)} \det(A)$
\item
$\det (A^s) = \det(A)^s$.
\end{itemize}

In general we will consider operators with nontrivial kernel.
So let $A$ be a densely defined operator and let $\ker(A^\infty) := \bigcup_{n\in \N} \ker(A^n)$ be the generalized kernel.
Let $P$ be the orthogonal projection onto $(\ker(A^\infty))^\perp$.
If the operator $A':= PA$ on the Hilbert space $(\ker(A^\infty))^\perp$ is admissible we will call $A$ admissible as well and define the {\bf reduced determinant} as
$$
\det'(A) := \det (A').
$$ 

{\bf Example:}
Consider the Lalpace operator $\lap := -(\der{x})^2$ on the sphere $S^1 \cong \R / 2\pi \Z$.
The Hilbert space $\CH$ will be $L^2(S^1)$ and $\CH'$ will be the span of $\{ e^{ikx} | k\in \Z\}$.
By the above it follows
$$
\det'(\lap ) = (2\pi)^2.
$$

Let $A$ be an admissible operator.
In the following we will be interested in the "characteristic function"
$\la \ra \det (A+\la)$ of $A$.

\begin{proposition}
Let $A$ be an admissible operator and let $\la\in\C$ with $\Re (\la)>0$.
Then the operator $A+\la$ is admissible again.
The function $\la \ra \det (A+\la)$ extends to an entire function with zeroes at the eigenvalues of $-A$, the order of a zero being its multiplicity as an eigenvalue.
For $\la \ra +\infty$ we have the following asymptotics (see \cite{Voros}):
$$
-\log \det (A+\la) \sim \sum_{\alpha_\nu \neq 0,-1,\dots} c_\nu \Ga (\alpha_\nu) \la^{-\alpha_\nu}
$$ $$
+\sum_{\alpha_\nu =-k\in\{ 0,-1,\dots\}}c_\nu \frac{(-1)^k}{k!} \left( \sum_{j=1}^k \rez{j} -\log \la\right) \la^k.
$$

Sometimes also the following is useful:
$$
-\log\det(A\pm i\la) = \sum_{\alpha_\nu=0} c_\nu (C_0+\ga \pm i\frac{\pi}{2} +\log\la)
$$ $$
+\sum_{\alpha_\nu =-k\in\{ 0,-1,\dots\}}c_\nu \frac{(\mp i)^k}{k!} \left( \sum_{j=1}^k \rez{j} -\log \la \mp i\frac{\pi}{2}\right)\la^k
$$ $$
+\sum_{\alpha_\nu \neq 0,-1,\dots} c_\nu \Ga (\alpha_\nu) (\pm i)^{-\alpha_\nu}\la^{-\alpha_\nu}+o(1),
$$
as $\la \ra +\infty$, where $C_0=\int_0^\infty \frac{1-\cos t}{t}dt$.

\end{proposition}

\prf
Let $M(s,\la)=\Ga (s)\zeta_{A+\la}(s)$. We have the differential equation:
$$
\der{\la}M(s,\la) = -M(s+1,\la).
$$
Now let $m\in \N$, we get
$$
(\der{\la})^{m+1} \zeta_{A+\la}(s) = (-1)^{m+1} (s+m)(s+m-1)\dots s\zeta_{A+\la}(s+m+1),
$$
so that for m large enough it follows
$$
(\der{\la})^{m+1} \zeta_{A+\la}(0) =0.
$$

It is clear that we have $\log \ \det (D+\la) = \lim_{s\rightarrow 0}(M(s,\la)-\zeta_{A+\la}(0))/s$ and the limit in $s$ may be interchanged with the $\la$-derivation. It follows that for $m$ large enough we have
$$
\begin{array}{cl} \displaystyle
(\der{\la})^{m+1}  \log \ \det(D+\la) & \displaystyle = (-1)^{m+1} M(m+1,\la)
\\
	& \displaystyle = (-1)^{m+1} \Ga(m+1)\sum_{n=1}^\infty \rez{(\la +\la_m)^{m+1}}.
\end{array}
$$
It remains to show the asymptotic expansion.
For this recall that we have $-\log \det (A+\la) = \frac{d}{ds}|_{s=0} \zeta_{A+\la}(s)$ and
$$
\zeta_{A+\la}(s) = \rez{\Ga (s)} \int_0^\infty t^{s-1} \Theta_{A+\la} (t) dt.
$$
As usual we split the integral as $\int_0^\infty = \int_0^1+\int_1^\infty$ and accordingly $\zeta_{A+\la}(s)=\zeta_{A+\la}^1(s)+\zeta_{A+\la}^2(s)$, where the second function is entire. Lets deal with $\zeta_{A+\la}^2(s)$ first. We estimate:
\begin{eqnarray*}
 |\frac{d}{ds}|_{s=0} \zeta_{A+\la}^2(s) | &=& | \int_1^\infty \Theta_{A+\la}(t) \frac{dt}{t}|\\
	&\leq & \sum_{n=0}^\infty \int_1^\infty e^{-t\Re (\la_n)}\frac{dt}{t}\ e^{-\la},
\end{eqnarray*}
which tends to zero exponentially and thus does not contribute to the asymptotics.
Now consider the integral $\int_0^1$.
The asymptotic expansion of $A$ tells us that for $N\in\N$ there is a $C>0$ such that
$$
|\Theta_A(t)-\sum_{\nu =0}^Nc_\nu t^{\alpha_\nu}| \leq Ct^{\alpha_{N+1}}
$$
for all $0<t<1$. This implies
$$
|\Theta_{A+\la}(t) - \sum_{\nu=0}^N c_\nu\sum_{n=0}^\infty \la^n \frac{(-1)^n}{n!}t^{\alpha_\nu +n}| \leq C t^{\alpha_{N+1}} e^{-\la t}.
$$
We can write
\begin{eqnarray*}
\zeta_{A+\la}^1(s) &=& \rez{\Ga (s)}\int_0^1 t^{s-1}\Theta_{A+\la}(t) dt\\
	&=& \sum_{\nu =0}^N \sum_{n=0}^\infty \la^n \frac{(-1)^n}{n!} c_\nu \rez{\Ga (s)} \rez{s+n+\alpha_\nu}\\
	&+& \rez{\Ga (s)}\int_0^1t^{s-1} \left( \Theta_{A+\la}(t) - \sum_{\nu=0}^N c_\nu \sum_{n=0}^\infty \la^n \frac{(-1)^n}{n!}t^{\alpha_\nu +n}\right) dt.
\end{eqnarray*}
For $N$ large enough this leads to
\begin{eqnarray*}
\frac{d}{ds}|_{s=0} \zeta_{A+\la}^1(s) &=& \sum_{\nu=0}^N \sum_{\begin{array}{c}n=0\\ n+\alpha_\nu=0\end{array}}^\infty \la^n \frac{(-1)^n}{n!} c_\nu \ga\\
	&+& \sum_{\nu=0}^N \sum_{\begin{array}{c}n=0\\ n+\alpha_\nu\neq 0\end{array}}^\infty \la^n \frac{(-1)^n}{n!} c_\nu \rez{n+\alpha_\nu}\\
	&+& \int_0^1 \left( \Theta_{A+\la}(t) - \sum_{\nu=0}^N c_\nu \sum_{n=0}^\infty \la^n \frac{(-1)^n}{n!}t^{\alpha_\nu +n}\right) \frac{dt}{t},
\end{eqnarray*}
where $\ga$ is the Euler constant.
We will treat the summands separately. The last summand is of absolute value less than
\begin{eqnarray*}
C \int_0^1 t^{\alpha_{N+1}} e^{-\la t} \frac{dt}{t} &=& C\la^{-\alpha_{N+1}}\int_0^\la t^{\alpha_{N+1}}e^{-t} \frac{dt}{t}\\
	&\leq & C\la^{-\alpha_{N+1}}\int_0^\infty t^{\alpha_{N+1}}e^{-t} \frac{dt}{t}\\
	&=& C\la^{-\alpha_{N+1}} \Ga (\alpha_{N+1}).
\end{eqnarray*}
Letting $N$ become large we see that the last summand does not contribute to the asymptotic expansion. Now the second summand splits into a contribution with $\alpha_\nu \in \{0,-1,-2,\dots\}$ and the complement. For $\alpha \in \C - \{ 0,-1,-2,\dots \}$ consider the function
$$
f_\alpha(x) := \sum_{n=0}^\infty \frac{(-x)^n}{n!} \rez{n+\alpha}.
$$
The auxilliary function $g_\alpha(x) :+\int_1^\infty t^{\alpha -1} e^{-xt} dt$ is rapidly decreasing in $x$ for any complex $\alpha$. A calculation shows that for $\alpha >0$ we have $f_\alpha(x) + g_\alpha(x) = x^{-\alpha}\Ga(\alpha)$ This gives meromorphic continuation to $\alpha \mapsto f_\alpha(x)$ and shows $f_\alpha \sim x^{-\alpha}\Ga(\alpha)$ for all $\alpha \in \C -\{ 0,-1,\dots \}$.

Now it remains to consider the terms with $\alpha_\nu =-k$. For these we consider
\begin{eqnarray*}
f_{-k}(x) &=& \sum_{k\neq n \geq 0} \frac{(-x)^n}{n!} \frac{1}{n-k}\\
	&=& \lim_{\alpha \ra -k}f_\alpha(x) - \frac{(-x)^k}{k!}\rez{k+\alpha}\\
	&=& -\log(x) x^k \frac{(-1)^k}{k!} + x^k a_k + g_{-k}(x),
\end{eqnarray*}
as we get by the de l'Hospital rule. Here $a_k$ is the derivative of $\Ga (s)(k+s)$ at the point $s=-k$.
This gives
$$
a_k = \frac{(-1)^{k+1}}{k!}\ga + \frac{(-1)^{k}}{k!} \sum_{j=1}^k \rez{k}.
$$
This gives the asymptotic expansion of $-\log\det(A+\la)$.
The formula for $-\log\det(A\pm i\la)$ follows similarly where we make extensive use of the fact that the Fourier transform of an $L^1$ function vanishes at infinity.
\qed

\subsection
There is also a more direct way to define determinants of operators on infinite dimensional spaces, namely the Fredholm determinant, which is defined as follows: Let T be a trace class operator, then we define
$$
\det_{Fr}(1+T) := \sum_{k=0}^\infty \tr \wedge^k T.
$$
The sum is absolutely convergent. Moreover, assume T normal with eigenvalues $(\lambda_n)_{n\in \N}$ then it follows that the infinite product
$$
\prod_{n\in \N} (1+\lambda_n)
$$
is convergent and equals the Fredholm determinant of $1+T$. If the operator norm $\parallel T \parallel$ of $T$ is less than 1 we also have
$$
\det_{Fr} (1+T) = \exp (-\sum_{n=1}^\infty \frac{(-1)^n}{n} \tr T^n).
$$

The connection between Fredholm determinant and zeta regularized determinant is:

\begin{proposition} Let $A$ and $A+1$ be admissible with $\ker (A)=0$. Assume $A^{\epsilon-1}$ is of trace class for some $\epsilon >0$, then
$$
\det_{Fr}(1+A^{-1}) = \frac{\det (A+1)}{\det (A)}.
$$
\end{proposition}

{\bf Proof:} Since the claim is certainly true for finite dimensions we can divide both sides by the factors stemming from the finitely many eigenvalues of A which are less than 1 in value. This means, we may assume $\parallel A^{-1} \parallel < 1$ so that we have at hand the exponential series description of the Fredholm determinant. Let $(\lambda_n)$ denote the eigenvalues of $T$ then we have:
$$
\begin{array}{cl} \displaystyle
\der{s} \zeta_{A+1} (s) & \displaystyle = -\sum_{n=1}^\infty \log(\la_n +1)(\la_n+1)^{-s}
\\ \displaystyle
	& \displaystyle = -\sum_{n=1}^\infty \log(\lambda_n) (\lambda_n +1)^{-s}
		- \sum_{n=1}^\infty \log (1+\rez{\lambda_n})(\la_n+1)^{-s}.
\end{array}
$$
The second term converges to $-\log \det_{Fr}(1+A^{-1})$ as $s\rightarrow 0$. The first equals
$$
-\sum_{n=1}^\infty \log (\la_n)\la_n^{-s} -\sum_{n=1}^\infty log(\la_n) \la_n^{-s} ((1+\rez{\la_n})^{-s}-1),
$$
in which the first summand is $\der{s}\zeta_A(s)$. Let $f_s(x)=(1+x)^{-s}-1$ then the Taylor series expansion of $f_s$ around $x=0$ is
$$
f_s(x) = \sum_{n=1}^\infty \frac{(-1)^n}{n!}s(s+1)\dots (s+n-1) x^n,
$$
so $f_s(x) =sxh(s,x)$ where $h$ is differentiable around $(x,s)=(0,0)$, so we have $|f_s(x)|\leq c_1sx$ for some $c>0$ thus
\begin{eqnarray*}
\left| \sum_{n=1}^\infty \log(\la_n) \la_n^{-s}((1+\rez{\la_n})^{-s}-1)\right|
&\leq & c_2s\sum_{n=1}^\infty \log(\la_n) \la_n^{-s-1}\\
&\leq & c_3 s \sum_{n=1}^\infty \la_n^{\epsilon -s-1}
\end{eqnarray*}
so this contribution vanishes at $s=0$.
\qed

\section{$L^2$-determinants}

The regularized determinant is defined via the zeta function $\zeta_A$ of an operator $A$ where $\zeta_A(s) = \tr A^{-s}$.
Suppose now, $A^{-s}$ belongs to some von Neumann algebra which is equipped with a canonical trace $\tau$.
 Then it would be natural to consider $\zeta_A^\tau(s) := \tau(A^{-s})$ as a generalization of the above.
If $\zeta_A^\tau$ again has a holomorphic continuation to zero then one gets the notion of the "$\tau$-determinant".
Seeking interesting trace functionals we will restrict our attention to an elliptic differential operator $D$ on a compact manifold.
We will however not perfectly stick to our philosophy since we will not find a  von Neumann algebra which contains the operator $D^{-s}$ but such an algebra will show up when we switch to the universal covering space and the lift of $D$ to this covering.

In order to have at hand the theory of heat kernels we will restrict to the following class of operators:
Let $E$ denote a smooth vector bundle over a smooth Riemannian manifold $X$.
A {\bf generalized Laplacian} on $E$ is a second order differential operator $D$ such that the principal symbol equals:
$$
\sigma_D (x,\xi) = \parallel \xi \parallel^2 \rm Id.
$$
We will further insist that $D$ is selfadjoint and semipositive, i.e. $D\geq 0$.
By \cite{BGV} it then follows that $D$ is admissible and thus $\det (D+\la)$ is well defined.

\subsection
From now on we fix the following situation: $X_\Ga$ denotes a smooth compact 
oriented Riemannian manifold with universal covering $X$ and fundamental 
group $\Ga$. 
We assume $\Ga$ infinite. Over $X_\Ga$ we will have a smooth Hermitian vector bundle $E_\Ga$ with pullback $E$ over $X$. On $E_\Ga$ we fix a generalized Laplacian $D_\Ga$ and we denote its pullback to $E$ by $D$.

Let $L^2(\Ga)$ denote the $L^2$-space over $\Ga$ with respect to the counting measure. The group $\Ga$ acts on $L^2(\Ga)$ by left translations. Inside the algebra $B(L^2(\Ga))$ of all bounded linear operators on the Hilbert space $L^2(\Ga)$ we consider the \index{von Neumann algebra of a group}{\bf von Neumann algebra of $\Ga$}, denoted $VN(\Ga)$ generated by the left translations $(L_\ga)_{\ga \in \Ga}$. There is a canonical trace on $VN(\Ga)$ defined by $\tr (\sum_\ga c_\ga L_\ga)=c_e$ which makes $VN(\Ga)$ a type ${\rm II}_1$ von Neumann algebra [GHJ], which is a factor if and only if every nontrivial conjugacy class in $\Ga$ is infinite.

It is easy to see that the commutant of $VN(\Ga)$ is the von Neumann algebra generated by the right translations.

Fix a fundamental domain ${\cal F}$ of the $\Ga$-action on $X$. We get isomorphisms of unitary $\Ga$-modules:
$$
L^2(E) \cong L^2(\Ga)\hat{\otimes}L^2(E{\mid_{\cal F}}) \cong L^2(\Ga) \hat{\otimes} L^2(E_\Ga).
$$
So for the von Neumann algebra $B(L^2(E))^\Ga$ of operators commuting with the $\Ga$-action we get
$$
B(L^2(E))^\Ga \cong VN(\Ga) \hat{\otimes}_w B(L^2(E_\Ga)),
$$
where $\hat{\otimes}_w$ means the weak closure of the algebraic tensor product, i.e. the tensor product in the category of von Neumann algebras.

Taking the canonical trace on $VN(\Ga)$ and the usual trace on $B(L^2(E_\Ga))$ we obtain a
type ${\rm II}_\infty$ trace on the algebra $B(L^2(E))^\Ga$ which we will denote by \index{$\tr_\Ga$} $\tr_\Ga$.
The corresponding dimension function is denoted by \index{$\dim_\Ga$} $\dim_\Ga$. Assume for example, a $\Ga$-invariant operator $T$ on $L^2(E)$ is given as integral operator with a smooth kernel $k_T$, then a computation shows
$$
\tr_\Ga(T) = \int_{\cal F} \tr(k_T(x,x))dx.
$$ \label{gamma-trace}

\subsection
Since the operator $D$ is the pullback of an operator on $\Ga \bs X$, it commutes with the $\Ga$-action and so does its heat operator $e^{-tD}$. This heat operator $e^{-tD} \in B(L^2(E))^\Ga$ has a smooth kernel $<x\mid e^{-tD}\mid y>$. We get:
$$
\tr_\Ga e^{-tD} = \int_{\cal F}\tr <x\mid e^{-tD}\mid x> dx.
$$
From this we read off that $\tr_\Ga e^{-tD}$ satisfies the same small time asymptotics as $\tr e^{-tD_\Ga}$ \cite{BGV}.

Let $D'=D\mid_{(ker(D))^\perp}$.
Unfortunately very little is known about large time asymptotics of $\tr_\Ga e^{-tD'}$ (see \cite{LL}). Let
$$
{\rm GNS}(D) := \sup \{ \alpha \in \R \mid \tr_\Ga e^{-tD'} = O(t^{-\alpha/2})\ {\rm as}\ t\rightarrow \infty \}
$$
denote the \index{Gromov-Novikov-Shubin invariant}{\bf Gromov-Novikov-Shubin invariant} of $D$ (\cite{GrSh}, \cite{LL}).
Then ${\rm GNS}(D)$ is always $\geq 0$. J. Lott showed in \cite{Lo} that the ${\rm GNS}$-invariants of Laplacians are homotopy invariants of a manifold. J.Lott and W. L\"uck conjecture in \cite{LL} that the Gromov-Novikov-Shubin invariants of Laplace operators are always positive rational or $\infty$.

\subsection
Throughout we will {\bf assume} that the Gromov-Novikov-Shubin invariant of $D$ is positive. We consider the integral
$$
\zeta_{D_\Ga}^1(s) := \rez{\Ga (s)} \int_0^1 t^{s-1} \tr_\Ga e^{-tD'}dt,
$$
which converges for $\Re (s) >>0$ and extends to a meromorphic function on the entire plane which is holomorphic at $s=0$, as is easily shown by using the small time asymptotics ([BGV],Thm 2.30).

Further the integral
$$
\zeta_{D_\Ga(s)}^2(s) := \rez{\Ga(s)}\int_1^\infty t^{s-1}\tr_\Ga e^{-tD'} dt
$$
converges for $\Re (s)<{\rm GNS}(D)$, so in this region we define the \index{$L^2$-zeta function}{\bf $L^2$-zeta function} of $D_\Ga$ as
$$
\zeta_{D_\Ga}^{(2)} (s) := \zeta_{D_\Ga}^1(s) + \zeta_{D_\Ga}^2(s).
$$

Assuming the Gromov-Novikov-Shubin invariant of $D$ to be positive we define the \index{$L^2$-determinant}{\bf $L^2$-determinant} of $D_\Ga$ as
$$
\det^{(2)}(D_\Ga) := \exp (-\frac{d}{d{s}} \mid_{s=0} \zeta_{D_\Ga}^{(2)} (s)).
$$

\begin{proposition}
Let $D_\Ga$ denote a generalized Laplacian over the manifold $X_\Ga$. Then the  \index{$L^2$-characteristic function}{\bf $L^2$-characteristic function}: $\la \mapsto \det^{(2)}(D_\Ga +\la),\ \la >0$, extends to a holomorphic function on $\C \backslash (-\infty ,0]$.
\end{proposition}

{\bf Proof:}
Let $\epsilon >0$ and consider the contour $\ga$ given by the negatively oriented boundary of the domain $G_\epsilon :=\{ \dist(z,]-\infty ,0])>\epsilon\}$.
For $\Re (s)>>0$, say $\Re (s)>R$ and $\la\in G_\epsilon$ we can write
$$
(D+\la)^{-s} = \frac{1}{2\pi i} \int_\ga z^{-s} (z -(D+\la))^{-1} dz .
$$
Varying $\epsilon$ this extends $(D+\la)^{-s}$ to a holomorphic pseudodifferential operator valued function on $\{\Re (s)>R\}\times \C -]-\infty ,0]$.
For $\Re (s) >R$ and $\Re (\la)>0$ we further have the differential equation
$$
\der{\la} \zeta_{D_\Ga +\la}^{(2)}(s) =-s \zeta_{D_\Ga +\la}^{(2)}(s+1).
$$
This integrates to
$$
\zeta_{D_\Ga +\la}^{(2)}(s) - \zeta_{D_\Ga +1}^{(2)}(s) = -s\int_1^\la \zeta_{D_\Ga +z}^{(2)}(s+1) dz.
$$
We thus extend the function $(\la ,s)\mapsto \zeta_{D_\Ga +\la}^{(2)}(s)$ to the domain $\C-]-\infty ,0] \times \{\Re (s)>R-1\}$.
Iterating this argument extends it to a meromorphic function on $\C-]-\infty ,0] \times \C$ which is holomorphic at $s=0$. This gives the claim.
\qed

\section{Convergence of determinants} 

In this section we will show that in case of a residually finite fundamental group  the characteristic functions of a tower will converge to the $L^2$-characteristic function.

\subsection
Let $\Ga$ be an infinite group. A \index{tower}{\bf tower} of subgroups is a sequence $\Ga_1 \supset \Ga_2 \supset \dots$ of subgroups of $\Ga$, each of which has finite index in $\Ga$ and the sequence satisfies $\cap_j \Ga_j = \{ 1\}$. A group $\Ga$ that has a tower is called \index{residually finite}{\bf residually finite}.

A tower $(\Ga_j)$ is called a {\bf normal tower} if each of the groups $\Ga_j$ is normal in $\Ga$.
Since the intersection of two finite index subgroups is a finite index subgroup it follows that to any tower $(\Ga_j)$ there exists a normal tower $\Ga_j'$ such that $\Ga_j'\subset \Ga_j$.
Since we think of $\Ga$ being a fundamental group we say that the normal tower $(\Ga_j')$ {\bf dominates} the tower $(\Ga_j)$.

Finitely generated subgroups of $GL_n(\C)$ are residually finite [Alp]. See [Kir] for more criteria of groups to be residually finite.

As a special example from arithmetic consider the following: Let $H$ denote a nonsplit quaternion algebra over $\Q$ which splits at $\infty$. Consider the group schemes $H^*$ and $H^1$ of units and norm one elements in $H$. Since via the splitting homomorphism $\pi : H^*(\R) \rightarrow GL_2(\R)$ the norm becomes the determinant, $H^1(\R)$ is mapped to $SL_2(\R)$. Let $R$ denote an order in $H$ then $\Ga (1) = \pi (H^1(R))$ is a cocompact discrete subgroup of $SL_2(\R)$. For an ideal $I$ in $R$ let
$ \Ga (I) = \{ \ga \in \Ga (1) \mid \ga \equiv 1\ \mod\ I \}$ the principal congruence subgroup. Now let $(I_n)_{n\in \N}$ be a decreasing sequence of nontrivial ideals in $R$ with $\cap_j I_j =0$ then the groups $\Ga_j =\Ga(I_j)$ form a tower.

\subsection
Now assume $\Ga = \pi_1(X_\Ga)$ is residually finite and fix a normal tower $(\Ga_j)_{j\in \N}$. For any continuous function $f$ on $\R$ we define the operator $f(D)$ by means of the spectral theorem as follows: Let $\mu$ denote a spectral resolution of $D$, i.e. $\mu$ is a projection valued measure on the spectrum ${\rm Spec} (D)$ of the operator $D$ such that on $L^2(E)$ we have
$$
D = \int_{{\rm Spec} (D)}\la d\mu (\la).
$$
We then define
$$
f(D) := \int _{{\rm Spec} (D)} f(\la) d\mu (\la).
$$
It follows
$$
\tr_\Ga(f(D)) = \int _{{\rm Spec} (D)} f(\la) d\tr_\Ga \mu (\la),
$$
where the integral will converge if $f(D)$ has a sufficiently smooth kernel.

Choose  $n\in \N$ larger then $2\dim G$, then $(1+D_{\Ga_j})^{-n}$ is of trace class and has a Schwartz kernel which is $(\dim X +1)$-times continuously differentiable. Let $F_{ev}^n$ denote the space of even $C^\infty$-functions on $\R$ 
which satisfy
$$
\parallel f\parallel_n ,\parallel f'\parallel_n ,
\parallel f''\parallel_n <\infty,
$$
where the norm is defined by $\parallel f\parallel_n :=\sup_{x\in\R}(1+x^2)^n|f(x)|$.

Note that the even Schwartz functions lie in $F_{ev}^n$.

\begin{theorem}
(Compare \cite{Don}.) For any $f\in F_{ev}^{n+1}$ we have
$$
\frac{\tr f(\sqrt{D_{\Ga_j}})}{[\Ga : \Ga_j]}
	\begin{array}{c}
		j\rightarrow \infty \\
		\longrightarrow		\\
		{}	
	\end{array}
	\tr_\Ga f(\sqrt{D}).
$$

For any Schwartz function $g\in \CS (\R)$ we have
$$
\frac{\tr g({D_{\Ga_j}})}{[\Ga : \Ga_j]}
	\begin{array}{c}
		j\rightarrow \infty \\
		\longrightarrow		\\
		{}	
	\end{array}
	\tr_\Ga g({D}).
$$
\end{theorem}

Before proving the theorem we give some immediate applications:

\begin{corollary}\label{cor1}
(Kazhdan inequality, compare \cite{Lu}) Let $h(D_{\Ga_j})$ denote the dimension of the kernel of $D_{\Ga_j}$ and $h^{(2)}(D_\Ga)$ the $\Ga$-dimension of $\ker (D)$ then we have
$$
\lim \ \sup_j \frac{h(D_{\Ga_j})}{[\Ga :\Ga_j]} \leq h^{(2)}(D_\Ga).
$$
\end{corollary}

We say that $D$ has a \index{spectral gap at zero}{\bf spectral gap at zero} if there is $\epsilon >0$ such that ${\rm Spec} D \cap ]0,\epsilon [=\emptyset$. Let
$$
N_j(x) := \sum_{\la \leq x} \dim {\rm Eig}  (D,\la ).
$$

\begin{corollary}\label{cor2}
Assume there is an $x>0$ such that
$$
\frac{N_j(x) -N_j(0)}
	{[\Ga :\Ga_j]} \jtoinfty 0.
$$then $D$ has a spectral gap at zero.
\end{corollary}

{\bf Proof of Corollary \ref{cor1}:}
Let $F$ denote the set of all $f\in \S (\R)$ with $f(0)=1$, $f(x)\geq 0$ for all $x\in \R$, then we have
$$
h(D_{\Ga_j})\leq \tr f(D_{\Ga_j})
$$
for all $f\in F$. It follows
$$
\limsup_j \frac{h(D_{\Ga_j})}{[\Ga :\Ga_j]}
	\leq \lim_j \frac{\tr f(D_{\Ga_j})}{[\Ga :\Ga_j]}
	= \tr_\Ga f(D)
$$
for all $f\in F$. So that
$$
\limsup_j \frac{h(D_{\Ga_j})}{[\Ga :\Ga_j]}
	\leq \inf_{f\in F} \tr_\Ga f(D) = h^{(2)}(D_\Ga).
$$\qed

{\bf Proof of Corollary \ref{cor2}:}
Assume $x>0$ as in the corollary. Let $f\in C^\infty_c(0,x)$ be positive, $f(y)\leq 1$ for all $y$, then
$$
N_j(x) -N_j(0) \geq \tr f(D_{\Ga_j}) \geq 0.
$$
Hence for all such $f$ we have $\tr_\Ga f(D) =0$ which gives the claim.
\qed

$ $

{\bf Proof of the theorem:}
Fix $x_0\in X$ and $v\in E_{x_0}^*$ the dual space to $E_{x_0}$. Recall that the distribution $u_s = \cos (s\sqrt{D}) v(\delta_{x_0}) \in C^\infty (E)'$
satisfies the wave equation
$$
(\frac{\partial^2}{\partial s^2} +D) u_s =0.
$$

By general results on hyperbolic equations (\cite{Tayl}, chap IV) it follows that $\cos (s\sqrt{D})$ has propagation speed $\leq \mid s \mid$. Now for $f\in F_{ev}^n$ the formula
$$
f(\sqrt{D}) = \rez{2\pi} \int_\R \hat{f} \cos (s\sqrt{D}) ds
$$
implies that $f(\sqrt{D})$ has finite propagation speed if its Fourier Transform $\hat{f}$ has compact support. Let $PW$ denote the Paley-Wiener space, i.e. the space of all Fourier transforms of smooth functions of compact support.
Then a standard verification shows that $PW$ is dense in $f\in F_{ev}^{n+1}$ with respect to the norm
$\parallel f \parallel_n$.

Let $f\in PW$, we will prove the assertion of the theorem for $f$.

Schwartz kernels of operators on sections of the bundle $E_\Ga$ are distributional sections of the bundle $E_\Ga \otimes_{ext} E_\Ga^*$ over $X_\Ga \times X_\Ga$. Where $\otimes_{ext}$ means the \index{exterior tensor product}{\bf exterior tensor product}, i.e.: $(E\otimes_{ext} F)_{(x,y)} \cong E_x \otimes F_y$.
Since $X_\Ga = \Ga \backslash X$, these can be identified with $\Ga \times \Ga$-invariant sections of $E\otimes_{ext} E^*$.

\begin{lemma}
Modulo the above identification we have for $f\in PW$ the identity of Schwartz kernels
$$
<\Ga x\mid f(\sqrt{D_\Ga})\mid \Ga y>
	= \sum_{\ga \in \Ga} <x \mid f(\sqrt{D})\mid \ga y>\ga .
$$
\end{lemma}

{\bf Proof of the lemma:}
Take $\ph \in C^\infty (E_\Ga)$ and identify $\ph$ with its pullback to $E$, which is a $\ga$-invariant section of $E$. The operator $\cos(s\sqrt{D})$ may be applied to $\ph$ in the distributional sense. Since $\cos (s\sqrt{D})$ is $\Ga$-invariant we obtain a $\Ga$-invariant distribution, so $\cos (s\sqrt{D})\ph \in C^\infty (E_\Ga)'$. we thus get an operator $\cos (s\sqrt{D})$ from $C^\infty (E_\Ga)$ to $C^\infty (E_\Ga)'$ satisfying the same differential equation and initial value conditions as $\cos (s\sqrt{D_\Ga})$, hence it equals the latter. For $f\in PW$ the formula $f(\sqrt{D}) = \rez{2\pi} \int_\R \hat{f} \cos (s\sqrt{D}) ds$ gives $f(\sqrt{D})\mid_{C^\infty(E_\Ga)} = f(\sqrt{D_\Ga})$. Since $f(\sqrt{D_\Ga})$ is a smoothing operator we may write
$$
\begin{array}{cl} 
f(\sqrt{D})\ph (x) &= \int_X <x\mid f(\sqrt{D})\mid y> \ph (y) dy
\\
	& = \sum_{\ga \in \Ga} \int_{\CF_\Ga} <x \mid f(\sqrt{D})\mid \ga y> \ga \ph (y) dy,
\end{array}
$$
where $\CF_\Ga$ denotes a fundamental domain of the $\Ga$-action on $X$. Now the sum
$$
\sum_{\ga \in \Ga} <x\mid f(\sqrt{D})\mid \ga y>\ga
$$
is locally finite and this gives the lemma.
\qed

To continue the proof of the theorem fix fundamental domains
$\CF_\Ga \subset \CF_{\Ga_1} \subset \CF_{\Ga_2} \subset \dots$ in a way that there are representatives for the classes in $\Ga /\Ga_j$ such that $\CF_{\Ga_j}=\cup_{\sigma : \Ga /\Ga_j}\sigma \CF_\Ga$.
We then get

\begin{eqnarray*}
\frac{\tr f(\sqrt{D_{\Ga_j}})}{[\Ga :\Ga_j]}
	&=& \rez{[\Ga :\Ga_j]} \sum_{\ga \in \Ga_j}
		\int_{\CF_{\Ga_j}} \tr <x \mid f(\sqrt{D}) \mid \ga x>\ga dx
\\
	&=& \sum_{\ga \in \Ga_j} \rez{[\Ga :\Ga_j]}
		\sum_{\sigma : \Ga /\Ga_j}
		\int_{\CF_{\Ga}} \tr
			<\sigma x \mid f(\sqrt{D}) \mid \ga \sigma x>\ga dx.
\end{eqnarray*}

The $\Ga$-invariance of the operator $f(\sqrt{D})$ implies $<\ga x |f(\sqrt{D})|\ga y> = \ga <x |f(\sqrt{D})| y>\ga^{-1}$ for all $\ga \in \Ga$. 
Thus the above equals
$$
\sum_{\ga \in \Ga_j} \rez{[\Ga :\Ga_j]}
		\sum_{\sigma : \Ga /\Ga_j}
		\int_{\CF_{\Ga}} \tr
			<x \mid f(\sqrt{D}) \mid \sigma^{-1}\ga \sigma x>\sigma^{-1} \ga\sigma dx.
$$
Since $\Ga_j$ is normal in $\Ga$ we end up with
$$
\frac{\tr f(\sqrt{D_{\Ga_j}})}{[\Ga :\Ga_j]}
	= \sum_{\ga \in \Ga_j} \int_{\CF_{\Ga}} \tr
			<x \mid f(\sqrt{D}) \mid \ga x>\ga dx.
$$
As $j\rightarrow \infty$ we get less and less summands, only the summand of $\ga =e$ remains, but this is just $\tr_\Ga f(\sqrt{D})$, so the claim follows for $f\in PW$.

Now let $f\in F_{ev}^n$ be arbitrary. 
Applying the above to $f(x)=e^{-tx^2}$ (see Lemma \ref{rapidly_decreasing}) and using standard techniques (\cite{D-stf}) it can be shown that there is a constant $C>0$ such that
$$
\frac{\tr(1+D_{\Ga_j})^{-n}}{[\Ga :\Ga_j]} \leq C \ \ \ {\rm for\ all}\ j\in \N.
$$

Let $\epsilon >0$. To our given $f\in F_{ev}^n$ there is a $g\in PW$ such that
$$
\sup_{x\in \R} \mid f(x) - g(x) \mid < \frac{\epsilon}{(1+x^2)^n}.
$$

We get

\begin{eqnarray*}
\mid \frac{\tr f(\sqrt{D_{\Ga_j}})}{[\Ga :\Ga_j]}
	- \tr_\Ga f(\sqrt{D}) \mid
	&\leq&  \mid \frac{\tr f(\sqrt{D_{\Ga_j}})}{[\Ga :\Ga_j]}
		- \frac{\tr g(\sqrt{D_{\Ga_j}})}{[\Ga :\Ga_j]} \mid
\\
	&+& \mid \frac{\tr g(\sqrt{D_{\Ga_j}})}{[\Ga :\Ga_j]}
		- \tr_\Ga g(\sqrt{D}) \mid
\\
	&+& \mid \tr_\Ga g(\sqrt{D}) - \tr_\Ga f(\sqrt{D}) \mid.
\end{eqnarray*}

The first summand on the right hand side is less than $\epsilon\frac{\tr((1+D_{\Ga_j})^{-n})}{[\Ga :\Ga_j]} \leq \epsilon C$. The second summand tends to zero as $j\rightarrow \infty$ and the third is less that $\epsilon \tr_\Ga (1+D)^{-n}$. The theorem follows.
\qed

For later use we prove:

\begin{lemma} \label{rapidly_decreasing}
The function
$$
t \mapsto \frac{\tr e^{-tD_{\Ga_j}}}{[\Ga :\Ga_j]}
	-\tr_\Ga e^{-tD},\ \ \ t>0,
$$
is rapidly decreasing at zero. Moreover this holds uniformly in j, that is, 
there is a function $f : ]0,1[ \rightarrow [0,\infty[$ which is rapidly 
decreasing at zero such that
$$
\mid \frac{\tr e^{-tD_{\Ga_j}}}{[\Ga :\Ga_j]}
	-\tr_\Ga e^{-tD} \mid \leq f(t)
$$
for all $j\in \N$ and all $t\in ]0,1[$.
\end{lemma}

{\bf Proof:}
Arguing as in the proof of the preceding lemma we see that
$$
<\Ga x\mid e^{-tD_\Ga}\mid \Ga y>
	= \sum_{\ga \in \Ga} <x\mid e^{-tD}\mid \ga y>\ga ,
$$
where convergence remains to be checked, however. To this end we write $d(x,y)$ for the Riemann distance of the points $x,y\in X$ and we use Theorem 3.1 of \cite{CGT} to get for $d(x,y)\geq a >0$:
$$
\mid <x\mid e^{-tD}\mid \ga y>\mid
	\leq \frac{C}{\sqrt{t}} e^{-(d(x,y)-a)^2/t},
$$
where $\mid . \mid$ here denotes the matrix $L^2$-norm, so $\mid A \mid = \sqrt{\tr AA^*}$.

By Proposition 4.1 in \cite{CGT} we conclude that the function 
$$
x\mapsto e^{-(d(x,y)-a)^2/t}
$$ 
is in $L^1(X)$ and the $L^1$-norm is clearly bounded by a constant independent of $y$. From this the convergence of the sum above follows. As in the previous section we get:
$$
\frac{\tr e^{-tD_{\Ga_j}}}{[\Ga :\Ga_j]} -\tr_\Ga e^{-tD_{\Ga}}
	= \sum_{1\neq \ga \in \Ga_j} \int_{\CF_\Ga}
		\tr<x \mid e^{-tD}\mid \ga x>\ga dx.
$$
Now fix $x_0\in \CF_\Ga$ and $r>0$ such that the ball of radius $r$ around $x_0$ contains $\CF_\Ga$, i.e. $\CF_\Ga \subset B_r(x_0)$. Choose $R>\max (r,1)$ such that $y\in \CF_\Ga ,x\in X$ with $d(x_0,x)>R$ implies $d(x,y)>r$.
Choose $j_0\in \N$ such that for all $j\geq j_0$ and all $\ga \in \Ga_j \backslash \{ 1\}$ we have $\bar{\ga \CF_\Ga} \cap B_R(x_0)=\emptyset$. Under these circumstances we get for $j\geq j_0$:

\begin{eqnarray*} 
\mid \frac{\tr e^{-tD_{\Ga_j}}}{[\Ga :\Ga_j]} -\tr_\Ga e^{-tD_{\Ga}} \mid
	&\leq & \sum_{\ga \in \Ga_j\backslash \{ 1\}} \frac{C}{\sqrt{t}}
		\int_{\CF_\Ga} e^{-(d(x,\ga x)-r)^2/t}dx
\\
	&\leq & \frac{C}{\sqrt{t}} \int_{X\backslash B_R(x_0)}
		e^{-(d(x,B_r(x_0))-r)^2/t}dx
\\
	&\leq & \frac{C}{\sqrt{t}} \int_R^\infty
		e^{-(\xi -r)^2/t} \vol (d(x,x_0) = \xi)d\xi
\\
	&\leq & \frac{C_1}{\sqrt{t}} \int_R^\infty
		e^{-(\xi -r)^2/t} e^{C_2 \xi} d\xi,
\end{eqnarray*}

for some $C_1, C_2 >0$ by Proposition 4.1 in \cite{CGT}.
But this can be calculated as

\begin{eqnarray*}
C_1 \int_{R/\sqrt{t}} e^{-(x-(\frac{r}{\sqrt{t}}+\frac{C_2 \sqrt{T}}{2}))^2
			+r +\frac{C_2^2 t}{4}} dx
&\leq& C_3 \int_{R/\sqrt{t}} e^{-x^2} dx
\\
&\leq & C_3 \int_{R/\sqrt{t}} x e^{-x^2} dx\ \ \ {\rm for}\ t<1
\\
&=& \frac{C_3}{2} e^{-R^2/t}
\end{eqnarray*}

which gives the assertion of the lemma.
\qed

For the following theorem we need not require the Gromov-Novikov-Shubin invariant to be positive.

\begin{theorem}
As $j\rightarrow \infty$ we have
$$
\det (D_{\Ga_j} +\la)^{\rez{[\Ga :\Ga_j]}}
	\longrightarrow \det^{(2)} (D_\Ga +\la)
$$
locally uniformly in $\la \in \C \backslash (-\infty ,0]$.
\end{theorem}

{\bf Proof:}
For $\Re \la >0$ and $s\in \C$ we consider
\begin{eqnarray*} 
F_j (\la ,s) &=& \frac{\zeta_{D_{\Ga_j}+\la}(s)}{[\Ga :\Ga_j]}
			- \zeta_{D_\Ga+\la}^{(2)}(s) \\
	&=& \rez{\Ga (s)} \int_0^\infty t^{s-1}
		(\frac{\tr e^{-tD_{\Ga_j}}}{[\Ga :\Ga_j]} -\tr_\Ga e^{-tD\Ga})
		e^{-t\la} dt
\end{eqnarray*}
which converges by the previous lemma. The formula
$$
F_j(\la ,s) = \frac{\tr (D_{\Ga_j}+\la)^{-s}}{\Ga :\Ga_j]}
		-\tr_\Ga (D+\la)^{-s}\ \ \ \Re (s) >>0
$$
together with the analytic continuation of the zeta function $\tr(D_{\Ga_j}+\la)^{-s}$ or $\tr_\Ga(D+\la)^{-s}$ 
shows that $F_j$ extends to a holomorphic function on the domain $(\C \backslash (-\infty ,0])\times \C$. As $\Re(\la)\rightarrow \infty$ we have $F_j(\la ,s)\rightarrow 0$ uniformly in $\{ \Re (s) < K\}$ for
all $K\in \R$. We want to show that this convergence is uniform in $j\in \N$. To this end let $\Re (\la) >0$ and split the integral above as $\int_0^1 +\int_1^\infty$. The integrand in the
$\int_0^1$-part obeys a uniform estimate by the previous lemma and thus the convergence is uniform in this part.

For the $\int_1^\infty$-part recall that there is a positive constant $C$ such that $\frac{\tr (D_{\Ga_j}+1)^{-n}}{[\Ga : \Ga_j]} < C$, from which it follows that
$\frac{\tr e^{-tD_{\Ga_j}}}{[\Ga :\Ga_j]}$ is bounded for $t\geq 1$ with a bound not depending on j. From this it follows
$$
F_j(\la ,s) \begin{array}{c} \Re (\la) \rightarrow \infty \\
		\longrightarrow \\ {} \end{array} 0
$$
locally uniformly in s and uniformly in $j\in \N$.

On the other hand we have a differential equation
$$
\der{\la} F_j(\la ,s) = F_j(\la ,s+1).
$$

For $\epsilon >0$ let $G_\epsilon$ denote the region
$$
G_\epsilon := \{ z\in \C \mid \Re (z)> -\rez{\epsilon} {\rm and}\
		{\rm dist}(z,(-\infty ,0])>\epsilon \}.
$$

From the above it follows that there is a $R>0$ such that for $\Re (s)>R$ we have
$$
F_j(\la ,s) \jtoinfty 0\ \ \ {\rm uniformly\ in}\ G_\epsilon
		\times \{ \Re (s) >R \} .
$$
We want to show that, having this convergence for $\Re (s)>R$ it already follows for $\Re (s)>R-1$. To achieve this let $\Re (s) > R-1$ and consider
$$
F_j(\la ,s)-F_j(\mu ,s) = \int_\mu ^\la F_j (z,s+1) dz,
$$
which tends to zero as $j\rightarrow \infty$. But on the other hand
$$
F_j (\la ,s) - F_j(\mu ,s) \begin{array}{c} \mu \rightarrow \infty \\
		\longrightarrow \\ {} \end{array} 	
		F_j(\la ,s),
$$
as the latter convergence is uniform in $j$ we get the claim. From this we get $F_j(\la ,s)\jtoinfty 0$ uniformly in a neighborhood of $s=0$ and this gives the theorem.
\qed

\subsection
Write $h(D_\Ga)$ for the dimension of the kernel of $D_\Ga$ and $h^{(2)}(D_\Ga)$ for the $\Ga$-dimension of the kernel of $D$. L\"uck showed in \cite{Lu}:
$$
\frac{h(\triangle_{p,\Ga_j})}{[\Ga :\Ga_j]} \jtoinfty h^{(2)}(\triangle_{p,\Ga})
$$
where $\triangle_p$ is the $p$-th Laplacian. Unfortunately the combinatorial arguments he used do not carry over to the analytic situation. So we only can show such a convergence under additional assumptions.

\begin{theorem}
Suppose that there is a number $\epsilon >0$ such that 
$$
{\rm Spec}(D_{\Ga_j})\cap ]0,\epsilon[ = \emptyset
$$ 
for all $j$, then
$$
\frac{h(D_{\Ga_j})}{[\Ga :\Ga_j]} \jtoinfty h^{(2)}(\Ga)
$$
and
$$
\det(D_{\Ga_j})^{\rez{[\Ga :\Ga_j]}} \jtoinfty \det^{(2)} (D_\Ga).
$$
\end{theorem}

Examples of this are Shimura manifolds $X_\Ga$, where $\Ga$ has Kazhdan property (T) \cite{Lub}.
This for example holds if the universal covering $X$ of $X_\Ga$ is a symmetric space with every simple factor of rank $>1$.

\prf
Apply the proof of the previous theorem to $D-\epsilon$. \qed

\begin{theorem}
Assume there exist constants $C,\alpha >0$ such that
$$
\frac{\tr e^{-tD_{\Ga_j}'}}{[\Ga :\Ga_j]} \leq C t^{-\alpha}
$$
for $t>1$ and all $j\in \N$, then the Gromov-Novikov-Shubin invariant of $D$ is $\geq 2\alpha$ and
$$
\frac{h(D_{\Ga_j})}{[\Ga :\Ga_j]} \jtoinfty h^{(2)}(\Ga)
$$
as well as
$$
\det (D_{\Ga_j})^\rez{[\Ga :\Ga_j]}\jtoinfty \det^{(2)} (D_\Ga).
$$
\end{theorem}

Examples for this are the Laplacians of flat tori or Heisenberg manifolds.

{\bf Proof:} We have that $\frac{\tr e^{-tD_{\Ga_j}}}{[\Ga :\Ga_j]}$ tends to $\frac{h(D_{\Ga_j})}{[\Ga :\Ga_j]}$ as $t\rightarrow \infty$ and by the assumption this is uniform in j. Thus from $\frac{\tr e^{-tD_{\Ga_j}}}{[\Ga :\Ga_j]}\jtoinfty \tr_\Ga e^{-tD}$ and $\tr_\Ga e^{-tD} \begin{array}{c} t\rightarrow \infty \\ \longrightarrow \\ {} \end{array} h^{(2)}(D_\Ga)$ the second assertion follows. The first is now easy and for the third recall that for $0<\Re (s) <\alpha$ we have
$$
\frac{\zeta_{D_{\Ga_j}}(s)}{[\Ga :\Ga_j]} - \zeta_{D_\Ga}^{(2)} (s)
	= \rez{\Ga (s)} \int_0^\infty t^{s-1}
		(\frac{\tr e^{-tD_{\Ga_j}}}{[\Ga :\Ga_j]} - \frac{h(D_{\Ga_j})}{[\Ga :\Ga_j]}
		-\tr_\Ga e^{-tD} +h^{(2)}(D_\Ga)) dt
$$

We split this integral into a $\int_0^1$- and a $\int_1^\infty$-part. The $\int_0^1$-part defines an entire function converging locally uniformly to zero in a neighborhood of $s=0$ as $j\rightarrow \infty$.
In the $\int_1^\infty$-part the integrand is dominated by some constant times $t^{-\alpha}$ which allows us to interchange the integration with the limit as $j\ra \infty$, so this part vanishes as $j\ra \infty$.
Together we get
$$
\frac{\zeta_{D_{\Ga_j}}(s)}{[\Ga :\Ga_j]} \jtoinfty \zeta_{D_\Ga}^{(2)} (s)
$$
uniformly around s=0 which gives the last claim.
\qed

$ $

\subsection
{\bf Question:}
Does for any generalized Laplacian the sequence 
$$
\frac{\dim \ker D_{\Ga_j}}{[\Ga :\Ga_j]}
$$ 
converge to $\dim_\Ga \ker D$?

\subsection
{\bf Question:}
Does for an arbitrary generalized Laplacian the sequence $\det (D_{\Ga_j})^\rez{[\Ga :\Ga_j]}$ converge to $\det^{(2)}(D_\Ga)$ if ${\rm GNS} (D_\Ga)>0$ ?

$ $

Let for $j\in \N$ and $x\in \R$:
$$
N_j (x) := \sum_{\la \leq x}\dim {\rm Eig} (D_{\Ga_j} ; \la)
$$
the \index{eigenvalue counting function}{\bf eigenvalue counting function} of $D_{\Ga_j}$. The condition of the last theorem would follow from the existence of constants $a,b,C>0$ such that
$$
\frac{N_j(x) -N_j(0)}{[\Ga :\Ga_j]} \leq C(x^a +x^b)\ \ \ x\geq 0,
$$
for all $j\in \N$. If this holds, we say the {\bf the counting functions satisfy a global growth estimate}.

In the following chapters we will specialize to the case where $X$ is a symmetric space without compact factors.

\subsection
{\bf Question:} Suppose $X$ is a symmetric space of the noncompact type and $D$ a generalized Laplacian which is invariant under the group of orientation preserving isometries of $X$. Do the counting functions of $D$ then satisfy a global growth estimate?

\section{A product formula}

\subsection
Suppose $g_{_\Ga} : X_\Ga \ra X_\Ga$ is an isometry.
Suppose further the action of $g_{_\Ga}$ lifts to a linear isometry of $E_\Ga$.
The {\bf $g_{_\Ga}$-equivariant zeta function} (see \cite{Koe} \cite{D-equiv})of $D_\Ga$ is defined by
$$
\zeta_{D_\Ga ,g_{_\Ga}}(s) := \sum_{\la >0} \tr (g_{_\Ga} | \Eig (D_\Ga ,\la)) \la^{-s}.
$$

\begin{lemma}
The equivariant zeta function extends meromorphically to the entire plane and  is holomorphic at $s=0$.
\end{lemma}

\prf
Use the fact that the operator $g_{_\Ga} e^{-tD_\Ga}$ has kernel
$$
<x | g_{_\Ga} e^{-tD_\Ga} | y> = <x | e^{-tD_\Ga} | g_{_\Ga} y>g_{_\Ga}
$$ 
and argue as in Lemma \ref{zetafunct}.
\qed

\subsection
We define the {\bf $g_{_\Ga}$-equivariant determinant} as
$$
\det_{g_{_\Ga}} (D_\Ga) := \exp (-\zeta_{D_\Ga ,g_{_\Ga}}'(0)).
$$

The following proposition is easy to see.

\begin{proposition}
The function 
$$
\la \mapsto \det_{g_{_\Ga}}(D_\Ga +\la),\ \ \ \la >0,
$$
 extends to a holomorphic function on $\C -(-\infty ,0]$.
We have $\det_{g_{_\Ga}}(D_\Ga) = \lim_{\la \downarrow 0}\det_{g_{_\Ga}}(D_\Ga +\la) \la^{-\tr (g_{_\Ga}|\ker D_\Ga)}$.
\end{proposition}
\qed

\subsection
Again we fix an isometry $g_{_\Ga}$ of $X_\Ga$ but now we also choose a lift $g$ to $X$.
This lift is an isometry of $X$ which is unique up to multiplication with elements of $\Ga$.
We are going to consider the operator $g e^{-tD}$.
The fact that the small time asymptotics hold pointwise \cite{BGV} implies that $\tr_\Ga (g e^{-tD})$ again satisfies a small time asymptotics.
The number $\tr_\Ga (g e^{-tD})$ does not depend upon the choice of the lift $g$.
Let
$$
GNS_g(D_\Ga) := \sup \{ \alpha\in\R | \tr_\Ga g e^{-tD'} = O(t^{\alpha /2})\}
$$
be the {\bf equivariant Gromov-Novikov-Shubin invariant}.
We will assume $GNS_g(D_\Ga)>0$.
Let
$$
\zeta_{D_\Ga ,g}^1(s) :+ \rez{\Ga (s)} \int_0^1 t^{s-1} \tr_\Ga (g e^{-tD}) dt,
$$
then this function is defined for $\Re (s) >>0$ and extends to a meromorphic function which is regular at $s=0$.
The integral
$$
\zeta_{D_\Ga ,g}^2(s) :+ \rez{\Ga (s)} \int_1^\infty t^{s-1} \tr_\Ga (g e^{-tD}) dt
$$
converges for $\Re (s) <\rez{2} GNS_g(D_\Ga)$. In this region we define the equivariant $L^2$-zeta function of $D_\Ga$ as
$$
\zeta_{D_\Ga ,g}^{(2)}(s) := \zeta_{D_\Ga ,g}^1(s) + \zeta_{D_\Ga ,g}^2(s),
$$
and the {\bf equivariant $L^2$-determinant} as
$$
\det^{(2)}_g(D_\Ga) :+ \exp (-\frac{d}{ds}|_{s=0} \zeta_{D_\Ga ,g}^{(2)}(s)).
$$

In analogy to the classical case we get

\begin{proposition}
The function $\la \mapsto \det_g^{(2)}(D_\Ga +\la)$, $\la >0$ extends to a holomorphic function on $\C -(-\infty ,0]$.
 
We have $\det_g^{(2)}(D_\Ga) = \lim_{\la \downarrow 0} \det_g^{(2)}(D_\Ga +\la)\la^{-\tr_{\Ga_g}(g|\ker D)}$.
\end{proposition}
\qed  

\subsection
Contrary to the non $L^2$-case the condition of $X_\Ga$ to be compact is not really needed in the $L^2$-setting.
One only needs the convergence of the integrals.
For $\ga\in\Ga$ let $\Ga_\ga$ denote the centralizer.
From \cite{D-equiv} we take

\begin{lemma} \label{tracegammagamma}
The operator $\ga e^{-tD}$ is of $\tr_{\Ga_\ga}$-trace class.
Its trace is
$$
\tr_{\Ga_\ga}(\ga e^{-tD}) = \int_{\CF_\ga}\tr <x | e^{-tD} | \ga x >\ga dx,
$$
where $\CF_\ga$ is a fundamental domain for $\Ga_\ga \bs X$ and the integral converges absolutely.

This integral can also be written as the integral over the compact set $X_\Ga$ of the smooth function
$$
x \mapsto \sum_{\tau \in [\ga]_\Ga}\tr <x|e^{-tD}|\tau x>\tau,
$$
where the sum runs over the $\Ga$-conjugacy class of $\ga$. 
Moreover, for $[\ga]\neq 1$ the function $t\mapsto \tr_{\Ga_\ga}(\ga e^{-tD})$ is rapidly decreasing for $t\downarrow 0$.
\end{lemma}
\qed

\begin{theorem}
For $\Re(\la) >0$ we have the locally uniformly convergent product expansion
$$
\det (D_\Ga+\la) = \prod_{[\ga]} \det_\ga^{(2)} (D_{\Ga_\ga}+\la),
$$
where the product runs over all conjugacy classes in $\Ga$.
\end{theorem}

\prf
As in the Proof of Lemma \ref{rapidly_decreasing} we have
$$
<\Ga x\mid e^{-tD_\Ga}\mid \Ga y>
	= \sum_{\ga \in \Ga} <x\mid e^{-tD}\mid \ga y>\ga 
$$
and therefore
\begin{eqnarray*}
\tr e^{-t(D_\Ga +\la)} &=& \sum_{\ga\in\Ga} \int_{\CF} \tr <x|e^{-t(D+\la)} |\ga x>\ga dx\\
	&=& \sum_{[\ga]} \sum_{\tau\in [\ga]} \int_\CF \tr <x|e^{-t(D+\la)} |\tau x>\tau dx\\
	&=& \sum_{[\ga]} \tr_{\Ga_\ga}(\ga e^{-t(D+\la)})
\end{eqnarray*}
according to Lemma \ref{tracegammagamma}.
By Lemma \ref{rapidly_decreasing} we conclude that for $\Re (\la)>0$ the difference of zeta functions
$$
\zeta_{D_\Ga +\la}(s) - \zeta_{D_\Ga +\la}^{(2)}(s)
$$
is entire and given by the integral
$$
\rez{\Ga (s)} \int_0^\infty t^{s-1} (\tr(e^{-t(D_\Ga +\la)}-\tr_\Ga e^{-t(D+\la)}) dt,
$$
which by the above is
$$
\sum_{[\ga]\neq 1} \zeta_{D_{\Ga_\ga} +\la}^{(2)}(s).
$$
The summands are equally well entire in $s$ and given by similar integrals.
Moreover the sum converges locally uniformly in $s$ and $\la$ which implies the claim.
\qed

\subsection
{\bf Example:}
Let $X_\Ga$ be a compact Riemannian Surface of genus $\geq 2$ with the hyperbolic metric  and let $\lap$ be the Laplace-Beltrami operator on $X_\Ga$.
The set of conjugacy classes $[\ga ]$ of $\Ga$ stands in a natural bijection to the set of all free homotopy classes of closed paths in $X_\Ga$.
Each such class has a member of minimal length $l_\ga$.
A closed path in the class of $[\ga ]$ is called {\bf primitive} if it is not a power of a shorter path.
This is a property of the class $[\ga ]$. 
To every $[\ga ]$ there is a unique primitive $[\ga_0 ]$ underlying $[\ga ]$.
We call the natural number $\mu_\ga := \frac{l_\ga}{l_{\ga_0}}$ the {\bf multiplicity} of $[\ga ]$.
By \cite{CaVo} we get for $[\ga ]\neq 1$:
$$
\det_\ga^{(2)} (D_{\Ga_\ga}+\la) = - \exp \left( \sum_{N\geq 0}\frac{e^{-(\sqrt{\la+\rez{4}}+N)}}{\mu_\ga}\right).
$$
For $[\ga ]=1$ \cite{CaVo} further shows
$$
\det_1^{(2)} (D_{\Ga_\ga}+\la)= \left( e^{-\la -\rez{4}} \det (P+\sqrt{\la +\rez{4}})\right)^{2-2g},
$$
where $P:=\sqrt{\lap^d +\rez{4}}$ and $\lap^d$ is the Laplace operator of the 2-sphere.

This shows that the determinant $\det (\lap +\la)$ equals
$$
\left( e^{-\la -\rez{4}} \det (P+\sqrt{\la +\rez{4}})\right)^{2-2g} \prod_{[\ga] prim} \prod_{N\geq 0} \left( 1-e^{-(\sqrt{\la +\rez{4}}+N)l_\ga}\right).
$$

\tiny
\today

\end{document}